\newcommand{\dataversione}{March 13, 2014}

\documentclass[11pt]{amsart}
\usepackage{mathrsfs}
\usepackage[notref,notcite]{showkeys}
\usepackage{hyperref}\hypersetup{colorlinks=true, citecolor=blue}
\usepackage{graphicx}
\usepackage{stmaryrd}
\numberwithin{equation}{section}

\newtheoremstyle{mytheorem}
{}
{}
{\it}
{\parindent}
{\bf}
{.}
{ }
{\thmnumber{#2.~}\thmname{#1}\thmnote{~\rm#3}}

\newtheoremstyle{myremark}
{}
{}
{\rm}
{\parindent}
{\bf}
{.}
{ }
{\thmnumber{#2.~}\thmname{#1}\thmnote{~\rm#3}}

\newtheoremstyle{myparagraph}
{}
{}
{\rm}
{\parindent}
{\bf}
{}
{ }
{\thmnumber{#2.~}\thmname{#1}\thmnote{#3}}

\theoremstyle{mytheorem}

\newtheorem{theorem}[subsection]{Theorem}

\theoremstyle{myremark}
\newtheorem{remark}[subsection]{Remark}
\newtheorem*{remark*}{Remark}

\theoremstyle{myparagraph}

\newtheorem*{parag*}{}

\makeatletter
\def\@secnumfont{\sc}
\def\section{\@startsection{section}{1}%
\z@{1.5\linespacing\@plus .2\linespacing}{.7\linespacing}%
{\normalfont\sc\centering}}
\makeatother

\makeatletter
\def\ps@headings{\ps@empty
 \def\@evenhead{%
  \setTrue{runhead}%
  \normalfont\footnotesize
  \rlap{\thepage}\hfil
  \def\thanks{\protect\thanks@warning}%
  \leftmark{}{}\hfil}%
 \def\@oddhead{%
  \setTrue{runhead}%
  \normalfont\footnotesize\hfil
  \def\thanks{\protect\thanks@warning}%
  \rightmark{}{}\hfil \llap{\thepage}}%
\let\@mkboth\markboth}
\makeatother
\makeatletter

\providecommand{\proofname}{Proof}
\makeatother
\newcommand{\R}{\mathbb{R}}

\newcommand{\dV}{d_V\kern-1pt}

\makeatletter
\@addtoreset{footnote}{section}
\makeatother

\begin{document}

	%
\pagestyle{empty}
\pagestyle{myheadings}
\markboth%
{\underline{\centerline{\hfill\footnotesize%
\textsc{Philippe Logaritsch and Andrea Marchese}%
\vphantom{,}\hfill}}}%
{\underline{\centerline{\hfill\footnotesize%
\textsc{Kirszbraun's extension theorem fails
for Almgren's multiple valued functions}%
\vphantom{,}\hfill}}}

	%
\thispagestyle{empty}

~\vskip -1.1 cm

	%
\centerline{\footnotesize version: \dataversione%
\hfill
}

\vspace{1.7 cm}

	%
{\Large\sl\centering
Kirszbraun's extension theorem fails\\
for Almgren's multiple valued functions\\
}

\vspace{.4 cm}

	%
\centerline{\sc Philippe Logaritsch and Andrea Marchese}

\vspace{.8 cm}

{\rightskip 1 cm
\leftskip 1 cm
\parindent 0 pt
\footnotesize

	%
{\sc Abstract.}
We show that there is no analog of Kirszbraun's extension theorem for Almgren's multiple valued functions.

\medskip\noindent
{\sc Keywords:} Kirszbraun's extension theorem, multiple valued functions, geometric measure theory.
\par
\medskip\noindent
{\sc MSC (2010): 54C20, 49Q20.}

\par
}

%
%
\section{Introduction}
Almgren's multiple valued functions play a key role in geometric measure theory since they are employed in the analysis of the branching behaviour of minimal surfaces in codimension larger than or equal to $2$ (see \cite{ALM} and \cite{DLS}).

We recall basic definitions for multiple valued functions. Let $Q$ be a positive integer, then
\[
 \mathcal{A}_Q(\R^n)=\left\{ \sum_{i=1}^Q\llbracket P_i\rrbracket : P_i\in\R^n,\, 1\leq i\leq Q \right\},
\]
where $\llbracket P\rrbracket$ denotes the Dirac measure at $P$. This space is endowed with the $L^2$-Wasserstein distance: for $T_1=\sum_{i=1}^Q\llbracket P_i \rrbracket$ and $T_2=\sum_{i=1}^Q\llbracket S_i\rrbracket$ we define
\[
 \mathcal{G}(T_1,T_2)=\min_{\sigma\in\mathcal{P}_Q}\sqrt{\sum_{i=1}^Q|P_i-S_{\sigma(i)}|^2},
\]
where $\mathcal{P}_Q$ denotes the group of permutations of $\{1,\dots, Q\}$.

One of the main ingredients in the theory of multiple valued funtions is the following extension theorem (see Theorem 1.7 in \cite{DLS}).
\begin{theorem}
\label{thm_DLS}
 Let $B\subset\R^m$ be a measurable set and let $f:B\to\mathcal{A}_Q(\R^n)$ be Lipschitz. Then there exists a constant $C=C(m,Q)>0$ and an extension $\bar{f}:\R^m\to\mathcal{A}_Q(\R^n)$ of $f$ such that
\[
\mathrm{Lip}(\bar{f})\leq C\mathrm{Lip}(f).
\]
\end{theorem}

In the Euclidean case, the classical Kirszbraun's extension theorem (see Theorem 2.10.43 in \cite{F}) states that an analogous result holds with $C=1$. 
More precisely, Kirszbraun's theorem states that Lipschitz functions defined on a subset of $\R^m$ with values in $\R^n$ (both endowed with the Euclidean distance) can be extended to all of $\R^m$ without increasing the Lipschitz constant.
The conclusion may fail as soon as $\R^m$ or $\R^n$ is remetrized by a metric which is not induced by an inner product, as shown in 2.10.44 in \cite{F}.

In \S\ref{cex} we prove that the conclusion also fails in the setting of multiple valued functions, by exhibiting a $\sqrt{2/3}$-Lipschitz function $f$ defined on a subset of $\R^2$ 
with values in $\mathcal{A}_2(\R^2)$ with the property that any Lipschitz extension $\bar{f}$ to $\R^2$ has Lipschitz constant at least 1. 

%
%

\section{Construction of the counterexample}
\label{cex}
Let $A=(0,1), B=(-\sqrt{3}/2,-1/2), C=(\sqrt{3}/2,-1/2)$ and let $P_1,\dots,P_6$ be the vertices of a regular hexagon centered at $0$, with side length 1: 
$P_1=(0,1)$, $P_2=(\sqrt{3}/2,1/2)$, $P_3=(\sqrt{3}/2,-1/2)$, $P_4=(0,1)$, $P_5=(-\sqrt{3}/2,-1/2)$ and $P_6=(-\sqrt{3}/2,1/2)$.

Consider the map $f:\{A,B,C\}\subset\R^2\to \mathcal{A}_2(\R^2)$ given by
\begin{eqnarray*}
 f(A)=&\llbracket P_1\rrbracket+\llbracket P_4\rrbracket,\\
f(B)=&\llbracket P_2\rrbracket+\llbracket P_5\rrbracket,\\
f(C)=&\llbracket P_3\rrbracket+\llbracket P_6\rrbracket.
\end{eqnarray*}
The Lipschitz constant of $f$ is $\sqrt{2/3}$. In fact, $|A-B|=|A-C|=|B-C|=\sqrt{3}$ and 
\[
\mathcal{G}(f(A),f(B))=\mathcal{G}(f(A),f(C))=\mathcal{G}(f(B),f(C))=\sqrt{2}.
\]
Now consider a map $\bar{f} :\{A,B,C\}\cup \{0\}\to\mathcal{A}_2(\R^2)$. 
We will prove that if $\bar{f}$ is an extension of $f$, then the Lipschitz constant of $\bar{f}$ is at least $1$.
Indeed, let $\bar{f}(0)=\llbracket S_1\rrbracket+\llbracket S_2\rrbracket$. Assume by contradiction $\mathrm{Lip}(\bar{f})< 1$, then $S_1$ and $S_2$ should lie on different sides of the perpendicular bisector of the line segment $\overline{P_1P_4}$ (see Figure 1).
In fact, if for example $S_1$ and $S_2$ both lie in the half plane $\{y\leq0\}$ then $|P_1-S_i|\geq 1$ for $i=1,2$ which implies $\mathcal{G}(\bar{f}(0),f(A))\geq 1$. 
The latter contradicts the assumption since $|A|=1$. 

\begin{figure}[ht]
 \includegraphics{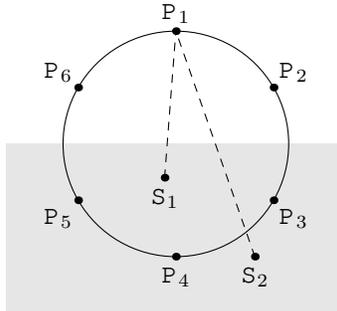}
\caption{$S_1$ and $S_2$ must lie on different sides of $y=0$}
\end{figure}

Arguing analogously for $\overline{P_2P_5}$ and $\overline{P_3P_6}$ we deduce that $S_1$ and $S_2$ must lie on opposite sectors among the six determined by the three perpendicular bisectors.
Without loss of generality we can assume that $S_1$ belongs to the intersection of the sector containing $P_1$ and the first orthant (see Figure 2).

\begin{figure}[ht]
 \includegraphics[scale=1]{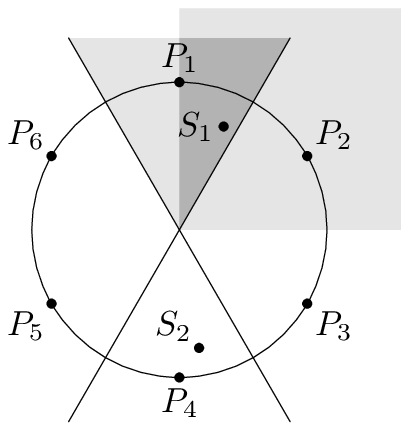}
\caption{}
\end{figure}

Since $|S_1-P_6|\leq|S_1-P_3|$ and $|S_2-P_6|\geq|S_2-P_3|$, we can estimate the distance between $\bar{f}(0)$ and $f(C)$ and get
\[
 \mathcal{G}(\bar{f}(0),f(C))^2=|S_1-P_6|^2+|S_2-P_3|^2\geq \frac{3}{4}+\frac{1}{4}=1,
\]
which contradicts our assumption since $|C|=1$.

\begin{remark}
Following the proof of Theorem \ref{thm_DLS} in \cite{DLS} one can explicitly determine the growth of the constant $C$ depending on $m$ and $Q$. 
It would be desirable to understand if the sharp constant has the same growth (or at least if $C(m,Q)$ goes to infinity as either $m$ or $Q$ goes to infinity).  
Clearly, just considering one-point extensions cannot lead to an answer to this question as the following general argument shows.

Let $(M,d_M)$ and $(N,d_N)$ be two complete metric spaces such that $M$ has the Heine-Borel property, $A$ a subset of $M$ and $f:A\to N$ be Lipschitz continuous. 
Then for every $P\in M\setminus A$ there exists a Lipschitz extension $\bar{f}:A\cup \{P\}\to N$ such that 
\[
 \mathrm{Lip}(\bar{f})\leq2 \mathrm{Lip}(f).
\]
In fact, let $S\in \bar{A}$ be a point realizing the distance between $P$ and $\bar{A}$. Let $\bar{f}(P)$ be the value at $S$ of the unique continuous extension of $f$ to $\bar{A}$, denoted by $f(S)$.
Then for every $y\in A\setminus \{S\}$ we get
\[
 \frac{d_N(\bar{f}(P),f(y))}{d_M(P,y)}=\frac{d_N(f(S),f(y))}{d_M(S,y)}\frac{d_M(S,y)}{d_M(P,y)}\leq 2 \mathrm{Lip}(f),
\]
because 
$d_M(S,y)\leq d_M(S,P)+d_M(P,y)$ and $d_M(S,P)\leq d_M(P,y)$ by the definition of $S$.

\end{remark}

%
%

\bibliographystyle{plain}

%
%

\vskip .5 cm

{\parindent = 0 pt\begin{footnotesize}

Ph.L. \& A.M. 
\\
Max-Planck-Institut f\"ur Mathematik
in den Naturwissenschaften
\\
Inselstrasse~22,
04103 Leipzig,
Germany
\\
e-mail Ph.L.: {\tt logaritsch@mis.mpg.de}\\
e-mail A.M.: {\tt marchese@mis.mpg.de}

\end{footnotesize}
}

\end{document}